\documentclass[10pt,a4paper]{article}

\usepackage{epsfig}
\usepackage{amsmath}
\usepackage{graphicx}
\usepackage{float}
\usepackage{subfig}
\usepackage{vmargin}
\usepackage{mathrsfs}
\usepackage{mathbbol}
\usepackage{hyperref}

\newcommand{\be}{\begin{equation}}
\newcommand{\ee}{\end{equation}}
\newcommand{\ba}{\begin{equation} \begin{aligned}}
\newcommand{\ea}{\end{aligned} \end{equation}}
\newcommand{\unitop}{\mathbb{1}}

\newcommand{\ket}[1]{\left| #1 \right>}

\author{Thomas House, Warwick Mathematics Institute} 
\title{Lie algebra solution of population models based on time-inhomogeneous
Markov chains} 
\date{}

\begin{document}

\maketitle

\begin{abstract}
Many natural populations are well modelled through time-inhomogeneous
stochastic processes. Such processes have been analysed in the physical
sciences using a method based on Lie algebras, but this methodology is not
widely used for models with ecological, medical and social applications. This
paper presents the Lie algebraic method, and applies it to three biologically
well motivated examples. The result of this is a solution form that is often
highly computationally advantageous.
\end{abstract}

\section{Introduction}

Stochastic models based on Markov chains are important in many ecological,
medical and social contexts. In these contexts, where populations are modelled,
there are often external influences that act on the system in a manner that
varies over time, leading to a time-inhomogeneous Markov chain model and
corresponding technical difficulties for analysis~\cite{Ross:2010}.

Wei and Norman~\cite{Wei:1963} proposed a method for dealing analytically with
time-inhomogeneous Markov chains, based on Lie algebraic methods. The idea of
combining Lie algebras and symmetry considerations with Markov chains has
continued to attract theoretical interest in a variety of
contexts~\cite{Wilcox:1967,johnson1985, mourad2004,Sumner:2011}.

At the same time, there is a more applied desire to have numerically efficient
methods to analyse Markov-chain population models, one option for which is the
use of matrix exponentials~\cite{Keeling:2008,Ross:2010}. The aim of this paper
is to explain how Lie algebraic methods can be used to derive matrix
exponential solutions to time-inhomogeneous Markov chains that are applicable
to population modelling.  In contrast to other applications, symmetries of
these systems are not a guide to the appropriate Lie algebra to use in solution
of population models; a certain amount of trial and error is necessary. The
focus of this paper is therefore on three examples in population modelling
where it is possible to define an appropriate Lie algebra, and a discussion of
the potential benefits of doing so.

\section{Methodology}

\subsection{Lie algebras}

In general, a Lie algebra over a field $F$ is an $F$-vector space $V$, together
with a bilinear map called the \textit{Lie bracket}. Elements of the vector
space are written $u,v,w\in V$. The Lie bracket is written $[u,v]\in V$ and
obeys
\be
[u,u] = 0 \text{ ,} \qquad [u,[v,w]] + [v,[w,u]] + [w,[u,v]] = 0 \text{ .}
\label{liedef}
\ee
We will be interested in the vector space $GL(n,{\mathbb{R}})$, i.e.\ the set
of real-valued $n\times n$ matrices, and will define the Lie bracket for
$X,Y\in GL(n,{\mathbb{R}})$ through the \textit{commutator}
\be
[X,Y] := XY - YX \text{ ,} \label{comdef}
\ee
which can be readily seen to satisfy~\eqref{liedef}. It is often convenient to
define an adjoint endomorphism operator, ad, to represent the Lie bracket:
\be
(\text{ad}X) Y := [X, Y] \text{ ,}
\label{addef}
\ee
so that multiple applications of the Lie bracket can be concisely written as
e.g.\ $(\text{ad}X)^2 Y = [X, [X,Y]]$ .

\subsection{Time-inhomogeneous Markov processes}

Suppose that $p(t)$ is a probability vector, i.e.\ a vector with values
$p_n(t)\geq 0$ summing to unity that represent the probability that an integer
stochastic variable takes the value $n$ at time $t$.  We consider models in
which the evolution of these probabilities over time is given by
\be
\frac{dp}{dt} = H(t) p(t) \text{ ,}
\label{KFeqn}
\ee
where $H(t)$ is a time-dependent matrix such that at any time $t$ its
off-diagonal elements are positive and its column sums are zero. This defines a
time-inhomogeneous continuous-time Markov chain. For some special cases,
analytic solutions can be obtained. But in general, where the state-space of
the Markov chain is finite, numerical algorithms exist that calculate $p(t)$ by
making use of expansions such as
\be
p(t + \delta t) = H(t)p(t) \delta t + O(\delta t^2) \text{ ,}
\label{euler}
\ee
and accumulating a sufficient number of $\delta t$ steps to reach $p(t)$ from
$p(0)$.  Methods based on this direct integration will therefore calculate
$p(t)$ in $O(t)$ operations.

\subsection{The method of Wei and Norman}

In this section, we recall the methodology of Wei and Norman~\cite{Wei:1963},
which is applicable to equations of the form~\eqref{KFeqn} above.  The first
step is to look for a decomposition of $H(t)$
\be
H(t) = \sum_{i=1}^{m} a_i(t) H_i \text{ ,}
\label{aexpand}
\ee
where the $H_i$ are linearly independent matrices obeying
\be
[H_i, H_j] = H_i H_j - H_j H_i = \sum_k {\xi_{ij}}^k H_k \text{ ,}
\label{liealg}
\ee
for (in our case real-valued) scalars ${\xi_{ij}}^k$.  Given such matrices, we
can form a vector space $V = \text{span}\{H_i\}_{i=1}^m \subseteq
GL(n,{\mathbb{R}})$ such that for a Lie bracket as defined in~\eqref{comdef},
we have $[X,Y]\in V$ for all $X,Y\in V$. This closure under the action of the
Lie bracket can be used to look for solutions of the form
\be
p(t) = e^{g_1(t) H_1} \cdots e^{g_m(t) H_m} p(0) =: U(t) p(0) \text{ ,}
\label{gexpand}
\ee
where matrix exponentiation is defined through the power series in the standard
way.  Using the ad operator as defined in~\eqref{addef}, the
Baker-Campbell-Hausdorff formula is
\be
e^{X} Y e^{-X} = e^{(\text{ad}X)} Y \text{ ,}
\label{BCH}
\ee
which will enable us to derive the solution form advertised.
Substituting~\eqref{aexpand} and~\eqref{gexpand} into~\eqref{KFeqn} then gives
\ba
\frac{dp}{dt} & = \sum_{i=1}^{m} a_i(t) H_i U(t) p(0) \\
& = \sum_{i=1}^{m} \dot{g}_i(t) \left( \prod_{j=1}^{i-1} e^{g_j (t) H_j} \right)
H_i \left( \prod_{j=i}^{m} e^{g_j (t) H_j} \right) p(0) \text{ .}
\ea
Since this expression holds for any $p(0)$, we can equate the operators acting on
$p(0)$, postmultiply by $U^{-1}$, and repeatedly apply~\eqref{BCH} to obtain
\be
\sum_{i=1}^{m} a_i(t) H_i = \sum_{i=1}^{m} \dot{g}_i(t)
 \left( \prod_{j=1}^{i-1} e^{g_j(t) (\text{ad}{H_j})} \right) H_i \text{ .}
\label{WN}
\ee
The precise solution to this equation will depend on the constants $\xi$
in~\eqref{liealg}, however since the $H_i$ are chosen to be linearly
independent, terms in~\eqref{WN} in front of the same basis matrix can be
equated, leading to a set of ODEs for the $g_i(t)$.

The usefulness of this method therefore depends on whether appropriate $H_i$
matrices can be defined, so that the equations that must be solved for $g_i(t)$
are not excessively complex. But in the event that $g_i(t)$ can be calculated
in $O(1)$ rather than, say, $O(t)$ -- which will typically be the case if an
analytic result is obtained -- then the computation of $p(t)$ can be achieved
in $O(1)$ rather than $O(t)$ through the numerical calculation of the matrix
exponentials in~\eqref{gexpand}. Such enhanced computational tractability of
stochastic models clearly has benefits in the application of probability theory
to statistical inference, where the speed of evaluation of likelihoods is
highly important.

\section{Examples}

The primary difficulty in applying the method above to population models is
finding an appropriate expansion of the form~\eqref{aexpand}, since the systems
involved are not obviously symmetric. We now turn to three examples where an
appropriate expansion can be found. In two cases, special initial conditions
give analytic results that can be checked against other methods to confirm the
soundness of the approach; and in the third, a significant numerical benefit is
observed compared to direct integration.

\subsection{A birth-death process}

Suppose we have a time-inhomogeneous birth-death process, characterised by a
stochastic variable $N(t)\geq 0$ taking integer values $n$, and the transition
rates
\ba
 n & \rightarrow n+1 \text{ at rate } b(t) \text{ ,}\\
 n & \rightarrow n-1 \text{ at rate } n d(t) \text{ .} \label{bdrate}
\ea
A biological interpretation of this process would be the survival of juvenile
animals, introduced to an inhospitable region by seasonal breeding happening at
another site, and dying at a rate that depends on the climate.  Defining
components of a vector $p_n(t) = \Pr(N(t)=n)$, the Kolmogorov equation for this
process is
\be
 \frac{dp}{dt} = (b(t) (R-\unitop{}) + d(t) (L - M)) p \text{ .}
\label{bdH}
\ee
The matrices involved are countably infinite in dimension, and are defined
implicitly by~\eqref{bdrate} and~\eqref{bdH}.  It is also possible to write
explicit definitions in terms of the Kronecker delta:
\be
 (\unitop{})_{n,k} = \delta_{n,k} \text{ ; } \ 
 (R)_{n,k} = \delta_{n,k+1} \text{ ; } \ 
 (L)_{n,k} = (k-1)\delta_{n,k-1} \text{ ; } \ 
 (M)_{n,k} = k\delta_{n,k} \text{ .}
\ee
Clearly, the identity matrix commutes with everything (i.e.\ $[\unitop{},X]=0$
for any $X$) while the other matrices obey
\be
[L,R] = \unitop{} \text{ ,}\quad
[M,R] = R \text{ ,}\quad
[L,M] = L \text{ .} \label{bdlie}
\ee
We then look for solutions of the form
\be
p(t) = e^{g_1(t) \unitop{}} e^{g_2(t) R} e^{g_3(t) L}
       e^{g_4(t) M} p(0) \text{ ,} \label{bdprod}
\ee
noting that $e^{g_1(t) \unitop{}} = e^{g_1(t)}$.  Making use of the
result~\eqref{WN} together with the algebra~\eqref{bdlie} gives
\ba
g_1(t) & = - e^{\mathcal{D}(t)} \int_0^t b(u) e^{\mathcal{D}(u)} du \text{ ,}\\
g_2(t) & = e^{\mathcal{D}(t)} \int_0^t b(u) e^{\mathcal{D}(u)} du \text{ ,}\\
g_3(t) & = e^{\mathcal{D}(t)} -1 \text{ ,}\\
g_4(t) & = -\mathcal{D}(t) \text{ ,}\\
\text{where}\quad \mathcal{D}(t) & := \int_0^t d(u) du \text{ .}
\label{bdgs}
\ea
This provides a solution to the original model, but one that is much simpler
if we assume the initial condition $N(0)=0$, in which case
\be
p_n(t) = e^{g_1(t)}\frac{\left(g_2(t)\right)^n}{n!} \text{ ,}
\label{bdsol}
\ee
for $g_1$, $g_2$ as in~\eqref{bdgs}. In this way an infinite-dimensional
time-inhomogeneous Markov chain is reduced to carrying out the two integrals
in~\eqref{bdgs}. It is worth comparing this to the `textbook' method for
dealing with time-inhomogeneous Markov chains, which is to derive an expression
for the probability generating function (PGF)~\cite{Grimmett:2001}.  This is
done by writing down the Kolmogorov equation~\eqref{bdH} in component form and
substituting in the definition of the PGF, $G(s,t) := \sum_n s^n p_n(t)$. This
gives a PDE for the PGF of
\be
\frac{\partial G}{\partial t} = (b(t)s - d(t))(s-1)
\frac{\partial G}{\partial s} \text{ ,}
\ee
which, given the initial condition $N(0)=0$, gives
\be
G(s,t) = \text{exp}\left[
(s-1) e^{-\mathcal{D}(t)} \int_0^t b(u) e^{\mathcal{D}(u)} du
\right] \text{ .} \label{bdpgf}
\ee
The equation~\eqref{bdpgf} yields the same solution as~\eqref{bdsol} above
through the standard relation
\be
p_n = \frac{1}{n!} \left. \frac{\partial^n G}{\partial s^n} \right|_{s=0}
\ee
The effort in deriving the solution~\eqref{bdsol} through the PGF and Lie
algebraic methods is therefore roughly similar; however the intermediate
results obtained in each method are likely to be useful in different contexts.
For example, the PGF in~\eqref{bdpgf} is likely to be the easiest way to derive
moments of the process; while the matrix exponential form~\eqref{bdprod} is
likely to be useful if the derivative of the solution with respect to a
parameter of the model is required~\cite{Wilcox:1967}.

\subsection{Epidemic surveillance}

Consider the following situation. An epidemic is in progress in a population,
such that individuals are either susceptible to infection, infectious, or
recovered and immune. Surveillance of the epidemic is carried out by
recruitment of individuals at random from the general population (or more
realistically through recruitment of individuals in contact with the healthcare
system due to non-infectious illness) who are tested and determined to be
either susceptible, infectious or recovered. The epidemic is characterised by a
force of infection $\lambda(t)$, which is the rate at which susceptible
individuals become infectious and for which there are various parametric
forms~\cite{Keeling:2007}, and also by a recovery rate $\gamma(t)$, which is
the rate at which infectious individuals recover. A plausible explicit choice
for these functions is to hold $\gamma(t)$ constant, and to take $\lambda(t)
=\lambda_0 e^{rt}$, representing the early exponential growth phase that is
common to many different epidemics.

As other authors have found, manipulation of more complex Markov chains is
simplified by the use of Dirac
notation~\cite{johnson1985,Jarvis:2005,Dodd:2009,Sumner:2011a}. In this
formalism, the probability vector for the model described above is written
\be
\ket{p(t)} = \sum_{S,I} \Pr (S,I|N,t) \ket{S,I} \text{ ,}
\ee
where $\Pr (S,I|N,t)$ is the probability that from a cohort of size $N$ a time
$t$ after the start of the epidemic a number $S$ of the cohort is susceptible
and a number $I$ is infectious (leaving $N-S-I$ recovered individuals).
$\ket{S,I}$ is a basis vector, linearly independent of any other basis vector
with different susceptible and infectious counts. Operators, marked out using a
hat $\hat{\mathcal{O}}$, act on basis vectors to give linear combinations of
basis vectors. For this system, we need the following operators:
\ba
\hat{S}\ket{S,I} & = S \ket{S,I} \text{ ,} \\
\hat{I}\ket{S,I} & = I \ket{S,I} \text{ ,} \\
\hat{\Delta}\ket{S,I} & = S \ket{S-1,I} \text{ ,} \\
\hat{\rho}\ket{S,I} & = I \ket{S,I-1} \text{ ,} \\
\hat{\tau}\ket{S,I} & = S \ket{S-1,I+1} \text{ .}
\ea
The action of these operators can be described computationally as follows:
$\hat{S}$ returns the number of susceptibles; $\hat{I}$ returns the number of
infectives; $\hat{\Delta}$ returns the number of susceptibles and depletes
these by one; $\hat{\rho}$ returns the number of infectives and depletes these
by one; and $\hat{\tau}$ returns the number of susceptibles, depletes the
susceptible population by one, and increases the infectious population by one.
The dynamical model is then
\be
\frac{d}{dt} \ket{p(t)} = \left( \gamma{}(t) \left(\hat{\rho} - \hat{I}\right)
 + \lambda{}(t) \left(\hat{\tau} - \hat{S} \right)\right) \ket{p(t)} \text{ .}
\label{epiH}
\ee
Note that while~\eqref{epiH} does not make use of the operator $\hat{\Delta}$,
it is necessary to include this to have an algebra that is closed under the
action of the Lie bracket. The full set of Lie brackets is shown in
Table~\ref{tab:sir}, and the action of the exponential adjoint endomorphism is
shown in Table~\ref{tab:sirad}. We then look for a solution of the form
\be
\ket{p(t)} = e^{g_1(t) \hat{\Delta}} e^{g_2(t) \hat{\tau}} e^{g_3(t) \hat{S}}
e^{g_4(t) \hat{\rho}} e^{g_5(t) \hat{I}} \ket{p(0)} \text{ .}
\label{sirexp}
\ee
Going through the same procedure as before gives solution
\ba
g_1(t) & = e^{\Lambda(t)} \left(1-e^{-\Lambda(t)}
- \int_0^t \lambda(u) e^{-\Lambda(u)} e^{\Gamma(u) - \Gamma(t)} du
\right)\text{ ,}\\
g_2(t) & = e^{\Lambda(t)}
\int_0^t \lambda(u) e^{-\Lambda(u)} e^{\Gamma(u) - \Gamma(t)} du
\text{ ,}\\
g_3(t) & = -\Lambda(t) \text{ ,}\quad
g_4(t) = 1-e^{-\Gamma(t)} \text{ ,}\quad
g_5(t) = -\Gamma(t)\text{ ,}\\
\text{for}\quad \Lambda(t) & := \int_0^t \lambda(u) du \text{ ,}\quad
\text{and}\quad \Gamma(t) := \int_0^t \gamma(u) du \text{ .}
\label{sirgs}
\ea
To check this result, assuming $\ket{p(0)}=\ket{N,0}$ and
substituting~\eqref{sirgs} into~\eqref{sirexp} gives
\be
\Pr (S,I|N,t) = \left( \begin{array}{c} 
N! \\ S! I! (N-S-I)!
\end{array} \right) (\pi_1)^S (\pi_2)^I (1-\pi_1-\pi_2)^{N-S-I}
\text{ ,}\label{pimulti}
\ee
where
\be
\pi_1 = e^{-\Lambda(t)} \text{ ,} \qquad
\pi_2 = \int_0^t \lambda(u) e^{-\Lambda(u)} e^{\Gamma(u) - \Gamma(t)} du
\text{ .}\label{pisol}
\ee
This is the form we would expect for this solution; considering each individual's
probability of remaining susceptible to be $\pi_1$ and being infections to be
$\pi_2$, these should obey
\be
\frac{d}{dt}\left( \begin{array}{c} \pi_1 \\ \pi_2 \end{array} \right) = 
\left( \begin{array}{cc} -\lambda(t) & 0 \\ \lambda(t)& - \gamma(t)
 \end{array} \right)\left( \begin{array}{c} \pi_1 \\ \pi_2 \end{array} \right)
\text{ ,} 
\ee
which has solution~\eqref{pisol}, and the independence of each individual leads
to the multinomial distribution~\eqref{pimulti} for the cohort as a whole. As
for the birth-death process, the analytic result obtained through Lie algebraic
methods can be obtained otherwise, and which method is preferable will depend on
which further calculations one wishes to undertake.

\subsection{A pure birth process}

\begin{figure}[H]
\begin{center}
\scalebox{0.9}{\resizebox{\textwidth}{!}{ \includegraphics{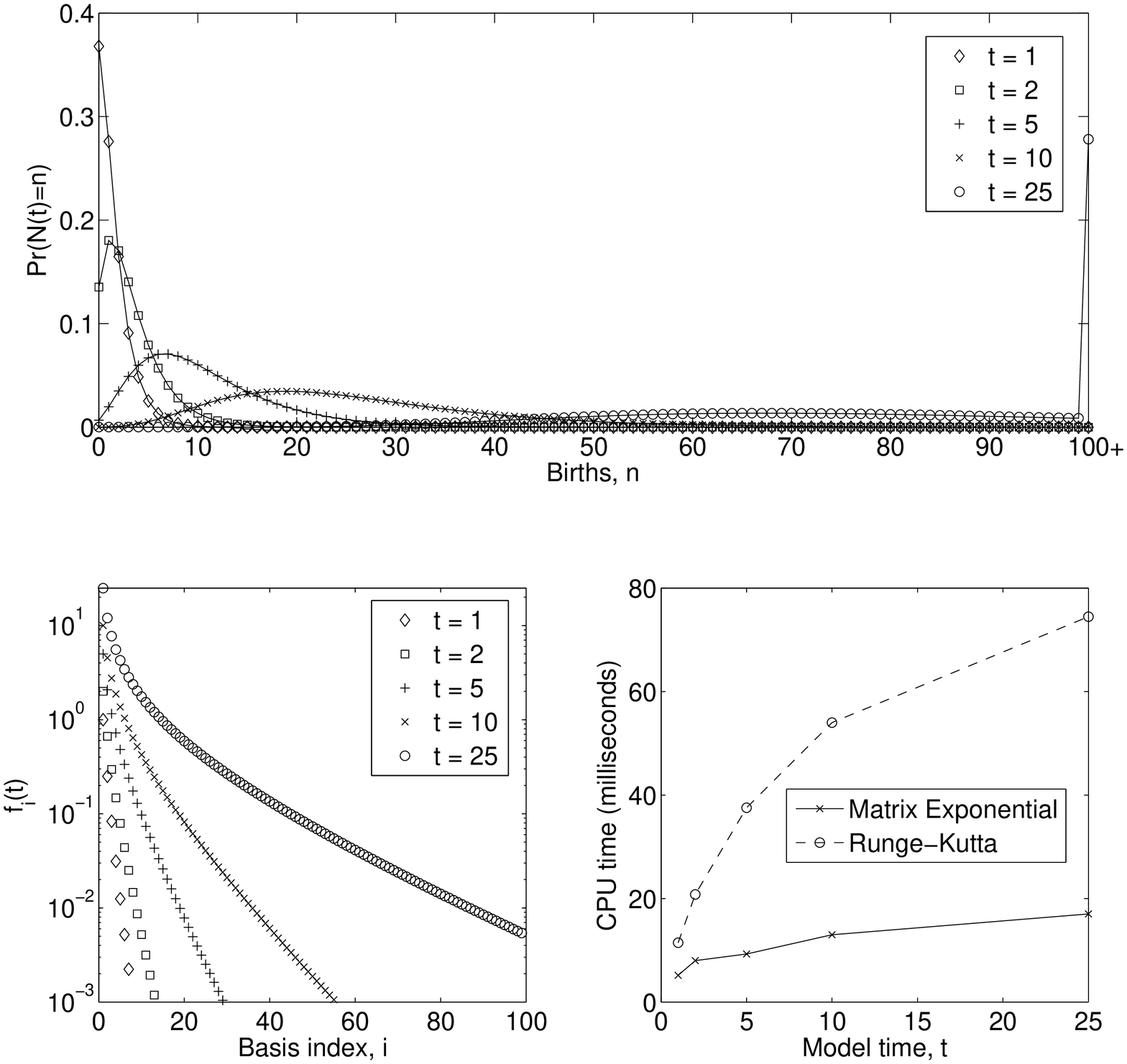} }}
\end{center}
\caption{Numerical results for a time-inhomogeneous pure birth process. Top:
probability distributions at different times. Matrix exponential solutions are
shown with markers and direct integration via Runge-Kutta is shown with
continuous black lines---clearly, these are numerically indistinguishable.
Bottom left: values of the quantities $f_i$ over time. Bottom right: CPU
time needed to run each method as a function of model time.}
\label{fig:birth}
\end{figure}

Now suppose we have a pure birth process, characterised by a stochastic
variable $N(t)$ and transition rate
\be
 n \rightarrow n+1 \text{ at rate } a(t) + n b(t) \text{ .} \label{pbrate}
\ee
Special forms of $a(t)$, $b(t)$ have been used to model the formation of social
contacts~\cite{danon:2011}. Defining components of a vector $p_n = \Pr(N(t)=n)$,
and assuming that there is a maximum count of interest $m$, such that we only
keep track of $\Pr(N(t)>m)$ as the final component of $p$, we write the
Kolmogorov equation for this process as
\be
 \frac{dp}{dt} = (a(t) P_1 + b(t) Q) p \text{ .}
\label{bdtimc}
\ee
The generating matrices $P_1$ and $Q$ are implicitly defined by~\eqref{pbrate},
however an additional set of matrices are required to produce a vector space
that is closed under the Lie bracket. These matrices take a similar form to
$P_1$ and are indexed by integer $i$. The definition of the matrices used, in
terms of the Kronecker delta, is
\be
 (Q)_{n,k} = k(\delta_{n,k+1} - \delta_{n,k}) \text{ ,} \qquad
 (P_i)_{n,k} = \delta_{n,k+i} - \delta_{n,k+i-1} \text{ .}
\ee
These matrices satisfy the commutation relations
\be
[P_i,P_j] = 0 \text{ ,} \qquad
[P_i,Q] = -iP_{i+1} + (i-1) P_i  \text{ .}
\label{pqcomm}
\ee
Note that it is here that the assumption of a finite state space (i.e.\ the
introduction of a maximum count of interest above) allows the solution method
to work, since otherwise a countably infinite number of matrices would be
needed.  We then look for a solution of the form
\be
p(t) = e^{f_1(t) P_1} \cdots e^{f_m(t) P_m} e^{g(t) Q} p(0) \text{ ,}
\ee
and make use of~\eqref{WN} to give the following system of equations:
\be
\dot{f}_1(t) = a(t) \text{ ,} \quad
\dot{g}(t) = b(t)  \text{ ,} \quad
\dot{f}_{i>1}(t) = (i-1)\dot{g}(t)\left[ f_{i-1}(t) - f_i(t) \right]\text{ .}
\ee
These equations have solution
\ba
{f}_1(t) & = \mathcal{A}(t) = \int_0^t a(u) du \text{ ,} \\
{g}(t) & = \mathcal{B}(t) = \int_0^t b(u) du   \text{ ,} \\
{f}_{i>1}(t) & = e^{-(i-1)\mathcal{B}(t)} \int_0^t a(u) 
 \left( e^{\mathcal{B}(t)} - e^{\mathcal{B}(u)}\right)^{i-1} du\text{ .}
\ea
Figure~\ref{fig:birth} shows numerical results for this system for the simple
choice
\be
 a(t) = 1 \text{ ,} \quad b(t) = \frac{1}{1+t} \text{ ,} \quad m = 100
\text{ ,} \quad N(0) = 0 \text{ .}
\ee
Note that for this initial condition, and due to the relations~\eqref{pqcomm},
we can write the solution as
\be
p_n(t)= \left( e^{\sum_i f_i(t) P_i} \right)_{n,0}
\text{ ,}
\ee
the evaluation of which which can be seen in Figure~\ref{fig:birth} to give a
significant computational advantage, as implemented in EXPOKIT~\cite{EXPOKIT},
compared to direct integration of~\eqref{bdtimc} through Runge-Kutta, as
implemented in the MATLAB function \texttt{ode45}, at large times.  Perhaps
unexpectedly, this is seen despite the relatively large value of $m$.  These
plots show that, while the ODE solver uses a more sophisticated relationship
than~\eqref{euler} to obtain better than $O(t)$ performance, it is still much
more sensitive to model time $t$ than the matrix exponential method.

There will, of course, be more complex systems where numerical integration via
Runge-Kutta is impractical, but analytic integration is simple (e.g.\ if either
$a(t)$ or $b(t)$ is a square wave rapidly oscillating between 0 and 1); for
these systems the matrix exponential method will further outperform
Runge-Kutta.  But there will also be systems where direct numerical integration
is straightforward and there is no simply obtained form for $f_i(t)$ and
$g(t)$, meaning that the matrix exponential solution is not useful.

\section{Discussion}

This paper has considered the solution of time-inhomogeneous Markov chains in
population modelling through the use of matrix exponentials. This is done using
the method of Lie algebras originally developed for applications in physical
sciences~\cite{Wei:1963}. In contrast to physical applications, population
models are often insufficiently symmetric to write down a well studied Lie algebra.
In even the relatively simple pure birth process considered, a large number of
basis matrices were needed to derive a matrix exponential solution; but despite
this the exponential solution is useful if a derivative with respect to a model
parameter is required~\cite{Wilcox:1967}, and perhaps more importantly often
has a significant numerical advantage over direct integration of the ODE
system~\cite{Keeling:2008}.  Given the popularity of computationally intensive
inference in modern population models~\cite{gilks:1995}, any such improvement
in numerical efficiency of likelihood evaluation can have important practical
benefits.

\section*{Acknowledgements}

Work funded by the UK Engineering and Physical Sciences Research Council. The
author would like to thank Josh Ross and Jeremy Sumner, in addition to the
editor and referee, for helpful comments on this work.

\newpage

\newpage

\begin{table}
\begin{center}
\begin{tabular}{cccccc}
\hline
$\hat{X}$ & $[\hat{X},\hat{S}]$ & $[\hat{X},\hat{I}]$ &
$[\hat{X},\hat{\Delta}]$ &
$[\hat{X},\hat{\rho}]$ & $[\hat{X},\hat{\tau}]$\\
\hline
$\hat{S}$ & 0 & 0 & $-\hat{\Delta}$ & 0 & $-\hat{\tau}$ \\
$\hat{I}$ & 0 & 0 & 0 & $-\hat{\rho}$ & $\hat{\tau}$ \\
$\hat{\Delta}$ & $\hat{\Delta}$ & 0 & 0 & 0 & 0 \\
$\hat{\rho}$ & 0 & $\hat{\rho}$ & 0 & 0 & $\hat{\Delta}$ \\
$\hat{\tau}$ & $\hat{\tau}$ & $-\hat{\tau}$ & 0 & $-\hat{\Delta}$ & 0 \\
\hline
\end{tabular}
\vspace{3em}
\caption{Values of $[\hat{X},\hat{Y}]$ for the epidemic model.}
\label{tab:sir}
\end{center}
\end{table}

\begin{table}
\begin{center}
\begin{tabular}{cccccc}
\hline
$\hat{X}$ & $e^{x(\text{ad}\hat{X})}\hat{S}$ &
$e^{x(\text{ad}\hat{X})}\hat{I}$ &
$e^{x(\text{ad}\hat{X})}\hat{\Delta}$ &
$e^{x(\text{ad}\hat{X})}\hat{\rho}$ & $e^{x(\text{ad}\hat{X})}\hat{\tau}$\\
\hline
$\hat{S}$ & $\hat{S}$ & $\hat{I}$ & $e^{-x}\hat{\Delta}$ & $\hat{\rho}$ &%
$e^{-x}\hat{\tau}$ \\
$\hat{I}$ & $\hat{S}$ & $\hat{I}$ & $\hat{\Delta}$ & $e^{-x}\hat{\rho}$ &%
$e^{x}\hat{\tau}$ \\
$\hat{\Delta}$ & $\hat{S} + x\hat{\Delta}$ & $\hat{I}$ & $\hat{\Delta}$ & $\hat{\rho}$ & $\hat{\tau}$ \\
$\hat{\rho}$ & $\hat{S}$ & $\hat{I}+x\hat{\rho}$ & $\hat{\Delta}$ & $\hat{\rho}$ & $\hat{\tau}+x\hat{\Delta}$ \\
$\hat{\tau}$ & $\hat{S}+x\hat{\tau}$ & $\hat{I}-x\hat{\tau}$ & $\hat{\Delta}$ & $\hat{\rho}-x\hat{\Delta}$ & $\hat{\tau}$ \\
\hline
\end{tabular}
\vspace{3em}
\caption{Values of $e^{x(\text{ad}\hat{X})}\hat{Y}$ for the epidemic model, for
scalar $x$.}
\label{tab:sirad}
\end{center}
\end{table}

\end{document}